\documentclass[12pt]{amsart}
\numberwithin{equation}{section}
\newtheorem{theorem}{Theorem}
\numberwithin{theorem}{section}

\numberwithin{theorem}{section} \numberwithin{lemma}{section}

\numberwithin{definition}{section}
\newtheorem{corollary}{Corollary}
\numberwithin{corollary}{section}

\numberwithin{remark}{section}

\numberwithin{proposition}{section}
\def\b{\begin{equation}}
\def\e{\end{equation}}
\newcommand{\ignore}[1]{}
\normalsize
\date {March 10, 2010}
\thanks{}
\begin{document}
\pagenumbering{arabic} \pagenumbering{arabic}\setcounter{page}{1}
\tracingpages 1
\title{Hardy-Poincar\'e, Rellich and Uncertainty principle  inequalities on Riemannian manifolds}
\author{Ismail Kombe and Murad \"Ozaydin}
\dedicatory {}
\address{Ismail Kombe, Mathematics Department\\ Dawson-Loeffler Science
\&Mathematics Bldg\\
Oklahoma City University \\
2501 N. Blackwelder, Oklahoma City, OK 73106-1493}
\email{ikombe@okcu.edu}
\address{Murad Ozaydin, Mathematics Department\\ University of
Oklahoma, Norman,  OK } \email{mozaydin@math.ou.edu}
\begin{abstract}
We continue our previous study of improved Hardy, Rellich and
Uncertainty principle inequalities on a Riemannian manifold $M$,
started in \cite{Kombe-Ozaydin}.
 In the present paper we prove new weighted Hardy-Poincar\'e,  Rellich type inequalities as well as improved version of our Uncertainty principle inequalities on  a Riemannian manifold $M$. In particular,  we obtain sharp constants for
 these inequalities on the hyperbolic space $\mathbb{H}^n$.
\end{abstract}
\maketitle
\section{Introduction}

The classical Hardy, Rellich and Heisenberg-Pauli-Weyl (uncertainty
principle) inequalities  play important roles in many questions
from spectral theory, harmonic analysis, partial differential equations, geometry
 as well as quantum mechanics. In order to motivate
our work, we present these three  classical (sharp) inequalities
on the Euclidean space $\mathbb{R}^n$. The Hardy inequality states
that for $n\ge 3$
\begin{equation}
\int_{\mathbb{R}^n}|\nabla\phi(x)|^2dx\ge
\Big(\frac{n-2}{2}\Big)^2\int_{\mathbb{R}^n}
\frac{|\phi(x)|^2}{|x|^2}dx,
\end{equation}
where $\phi\in C_0^{\infty}( \mathbb{R}^n)$.  Here the constant
$(\frac{n-2}{2})^2$ is  sharp, in the sense that

\[\begin{aligned}\Big(\frac{n-2}{2}\Big)^2 &=
\inf_{0\neq \phi\in
C_0^{\infty}(\mathbb{R}^n)}\frac{\int_{\mathbb{R}^n}|\nabla
\phi(x)|^2dx}{\int_{\mathbb{R}^n}\frac{|\phi(x)|^2}{|x|^2}dx}.
\end{aligned}\]

Another inequality involving  second order derivatives is the
Rellich inequality \cite{Rellich}:
\begin{equation} \int_{\mathbb{R}^n}|\Delta \phi(x)|^2dx\ge
\frac{n^2(n-4)^2}{16}\int_{\mathbb{R}^n
}\frac{|\phi(x)|^2}{|x|^4}dx,\end{equation} where $\phi\in
C_0^{\infty}(\mathbb{R}^n$), $n \ge 5$ and the constant
$\frac{n^2(n-4)^2}{16}$   is again sharp.( There are also versions
for lower dimensions under additional hypotheses.)

The classical Heisenberg-Pauli-Weyl  inequality, a precise mathematical
formulation of the uncertainty principle of quantum mechanics,
states that:
\begin{equation}
\Big(\int_{\mathbb{R}^n} |x|^2
|f(x)|^2dx\Big)\Big(\int_{\mathbb{R}^n} |\nabla f(x)|^2 dx\Big)\ge
\frac{n^2}{4} \Big(\int_{\mathbb{R}^n} |f(x)|^2 dx \Big)^2
\end{equation}
for all $f\in L^2( \mathbb{R}^n)$. Here the constant $\frac{n^2}{4}$
is sharp and  also it is well-known that equality is attained in
(1.3) if and only if $f$ is a Gaussian (i.e.
$f(x)=Ae^{-\alpha|x|^2}$ for some $A \in \mathbb{R}, \alpha>0$).

These inequalities  have been extensively studied in the Euclidean
setting  and now the literature on this topic is quite vast and
rich, encompassing many generalizations and refinements, e.g. \cite{Allegretto}, 
\cite{Folland-Sitaram}, \cite{Brezis}, \cite{Adimurthi},
\cite{Barbatis-Filippas-Tertikas}, \cite{Moradifam-1}, \cite{Cowan},
\cite{Moradifam-2}
 and references therein. Many new
developments are still forthcoming. For instance, Tertikas and Zographopoulos \cite{Tertikas} give a  sharp
Rellich-type inequality and its improved versions
which involves both first and  second order derivatives:
\begin{equation}
\int_{\mathbb{R}^n}|\Delta \phi(x)|^2dx\ge
\frac{n^2}{4}\int_{\mathbb{R}^n}\frac{|\nabla \phi(x)|^2}{|x|^2}dx,
\end{equation}
where $\phi\in C_0^{\infty}(\mathbb{R}^n)$, $n\ge 5$ and  the
constant $\frac{n^2}{4}$ is sharp.

On the other hand  the Euclidean results mentioned above continues
to be a source of inspiration for the problem of finding analogues
 inequalities in the setting of Riemannian manifolds. There
has been continuously growing literature in this direction, e.g.
\cite{Carron}, \cite{Davies-Hinz}, \cite{Grillo}, \cite{Barbatis},
\cite{Li-Wang}, \cite{Xia}, \cite{Kombe-Ozaydin}, \cite{Minerbe},
and the references therein. For instance, in an interesting paper
Carron \cite{Carron} obtained the following weighted $L^2$-Hardy
inequality on a complete non-compact Riemannian manifold $M$:

\begin{equation}\int_M \rho^{\alpha}|\nabla \phi|^2dV\ge
\Big(\frac{C+\alpha-1}{2}\Big)^2\int_
M\rho^{\alpha}\frac{\phi^2}{\rho^2}dV\end{equation} where
 $\phi\in
C_c^{\infty}(M-\rho^{-1}\{0\})$, $\alpha\in \mathbb{R}$, $C>1$,
$C+\alpha-1>0$ and the weight function $\rho$ satisfies  $|\nabla
\rho|=1$ and $\Delta\rho\ge \frac{C}{\rho}$ in the sense of
distribution. For complete non-compact Riemannian manifolds, under
the same geometric assumptions on the weight function $\rho$ we
obtained in \cite{Kombe-Ozaydin} an $L^p$ version of (1.5) (where
$ 1< p<\infty$ and $ C+1+\alpha-p>0$):
\begin{equation} \int _M\rho ^{\alpha}|\nabla\phi|^p dV\ge
\Big(\frac{C+1+\alpha-p}{p}\Big)^p\int _M
\rho^{\alpha}\frac{|\phi|^p}{\rho^ p}dV,\end{equation}  as well as
a Rellich-type inequality (where $\alpha<2$, $C+\alpha-3>0$):

\begin{equation}\int_{M}\rho^{\alpha}|\Delta \phi|^2dV \ge
\frac{(C+\alpha-3)^2(C-\alpha+1)^2}{16}
\int_{M}\rho^{\alpha}\frac{\phi^2}{\rho^4}dV\end{equation}
where $\Delta$ is the Laplace-Beltrami operator on $M$.
\medskip

We also found an $L^p$ Heisenberg-Pauli-Weyl uncertainty principle
type inequality (for a complete noncompact Riemannian manifold) and
an $L^2$ version with a (nonnegative) remainder term. In the
specific case when the manifold $M$ is the hyperbolic space
$\mathbb{H}^n$, we obtained sharp constants  for the Hardy and
Rellich-type inequalities, and  explicit (not sharp) constants for
the Heisenberg-Pauli-Weyl uncertainty inequalities.
\medskip

In the present paper we continue our investigation on Hardy, Rellich
and Heisenberg-Pauli-Weyl type inequalities. The plan of the paper
is as follows. In Section 2 we first prove a new form of weighted
Hardy-Poincar\'e type inequality and then  we prove various improved
versions of the weighted Hardy inequality (1.5)( in the sense that
nonnegative terms are added in the right hand side of (1.5)). We
note that these improved inequalities are the main  tool in proving
improved Rellich type inequalities. In Section 3 we first prove a
weighted analogue of (1.4) and then obtain improved versions.
Section 4 is devoted to the study of Heisenberg-Pauli-Weyl
(uncertainty principle) type inequalities where we obtain  better
constants than those  of \cite{Kombe-Ozaydin} and prove sharp
analogue of the classical  uncertainty principle inequality (1.3) on
the Hyperbolic space $\mathbb{H}^n$. In each section we first prove
inequalities in the context of a general complete Riemannian
manifold. Then, turning our attention to hyperbolic space
$\mathbb{H}^n$,  we consider specific weight functions and obtain
inequalities with explicit and usually sharp constants.

\section{Weighted Hardy-Poincar\'e type inequalities}

Throughout this paper, $M$ denotes a complete noncompact Riemannian manifold  endowed with a metric $g$. We denote by $dV$, $\nabla$, and $\Delta$ respectively the Riemannian volume element, the Riemannian gradient and the Laplace-Beltrami operator on $M$.

We begin this section by proving a new form of the Hardy-Poincar\'e
type inequality for a complete noncompact Riemannian manifold $M$
with a weight function $\rho$ modelled on the distance from a point.
(In this context the hypotheses $|\nabla \rho|=1$ and $\Delta\rho\ge
\frac{C}{\rho}$ seem to be geometrically quite natural.) One
advantage of this set-up is that it implies and thus provides
another (shorter) proof of (1.6) above (Theorem 2.1 in
\cite{Kombe-Ozaydin}) as explained in the \textbf{Remark} below.

\begin{theorem} Let $M$ be a complete noncompact Riemannian manifold of dimension
$n>1$. Let $\rho$ be a nonnegative function on $M$ such that
$|\nabla\rho|=1$ and $\Delta\rho\ge \frac{C}{\rho}$ in the sense of
distribution where $ C>0$. Then the following inequality hold:

\begin{equation}\int_{M} \rho^{\alpha+p}|\nabla \rho\cdot \nabla\phi|^pdx
\ge \Big(\frac{C+\alpha+1}{p}\Big)^p \int_{M}
\rho^{\alpha}|\phi|^pdx.\end{equation} for all compactly supported
smooth functions $\phi\in C_0^{\infty}(M\setminus
\rho^{-1}\{0\})$, $1<p<\infty$, and $C+\alpha>-1$.
\end{theorem}
\proof

It follows from above hypothesis that \begin{equation}
\text{div}(\rho\nabla \rho)\ge C+1.\end{equation} Multiplying both
side
 of (2.2) by $\rho^{\alpha}|\phi|^p$ and integrating over $M$ yields
\[(C+1)\int_M \rho^{\alpha}|\phi|^p dx\le \int_M \text{div}(\rho\nabla
\rho)\rho^{\alpha}|\phi|^pdx.\] As an immediate consequence of
divergence theorem we have

\[(C+\alpha+1)\int_M\rho^{\alpha}|\phi|^pdx \le -p\int_M |\phi|^{p-2}\phi \rho^{\alpha+1}
\nabla \rho\cdot\nabla \phi dx.\] An application of H\"older's and
Young's inequality yields
\[\begin{aligned}(C+\alpha+1)\int_M \rho^{\alpha}|\phi|^pdx &\le p
\Big(\int_M \rho^{\alpha}|\phi|^pdx\Big)^{(p-1)/p}\Big(\int_M
\rho^{\alpha+p}|\nabla\rho\cdot
\nabla\phi|^p dx\Big)^{1/p}\\
&\le (p-1)\epsilon^{-p/(p-1)}\int_M \rho^{\alpha}|\phi|^pdx+
\epsilon^p \int _M \rho^{\alpha+p}|\nabla \rho\cdot \nabla\phi|^p dx
\end{aligned}\] for any $\epsilon>0$. Therefore

\begin{equation}
\int_M \rho^{\alpha+p}|\nabla \rho\cdot \nabla \phi|^p dx\ge
\epsilon^{-p}\Big(C+\alpha+1-(p-1)\epsilon^{-p/(p-1)}\Big) \int_M
\rho^{\alpha}|\phi|^pdx.
\end{equation} Note that the function $\epsilon \longrightarrow
\epsilon^{-p}\Big(C+\alpha+1-(p-1)\epsilon^{-p/(p-1)}\Big)$
attains the maximum for $\epsilon^{p/(p-1)}=\frac{p}{C+\alpha+1}$,
and this maximum is equal to  $\Big(\frac{C+\alpha+1}{p}\Big)^p$.
Now we obtain the desired inequality: \[\int_M
\rho^{\alpha+p}|\nabla \rho\cdot \nabla \phi|^p dx \ge
\Big(\frac{C+\alpha+1}{p}\Big)^p \int_M \rho^{\alpha}|\phi|^pdx.\]
\qed

 \noindent\textbf{Remark}. Applying the Cauchy-Schwarz
inequality to $|\nabla\rho\cdot \nabla\phi|$,  replacing $\alpha$
with $\alpha-p$ and using $|\nabla\rho|=1$ yields the weighted
$L^p$-Hardy inequality (1.6).
\medskip

We will give a sharp version of Theorem 2.1 in  the hyperbolic space
$\mathbb{H}^n$. Recall that the hyperbolic space $\mathbb{H}^n$
$(n\ge 2)$ is a complete simple connected Riemannian manifold having
constant sectional curvature equal to $-1$. There are several models
for $\mathbb{H}^n$ and we will use the Poincar\'e ball model
$\mathbb{B}^n$ in this paper.
\medskip

The Poincar\'e ball model for the hyperbolic space is:
\[\mathbb{B}^n=\{x=(x_1, \cdots, x_n) \in \mathbb{R}^n |\, |x|<1\}\]
endowed  with the Riemannian metric $ds=\lambda(x)|dx|$ where
$\lambda(x)=\frac{2}{1-|x|^2}$. Hence $\{\lambda dx_i\}_{i=1}^{n}$
give an orthonormal basis of the tangent space at $x=(x_1,\cdots,
x_n)$ in $\mathbb{B}^n$. The corresponding dual basis is
$\{\frac{1}{\lambda}\frac{\partial}{\partial x_i}\}_{i=1}^{n}$,
thus the hyperbolic gradient and the Laplace Beltrami operator
are:
\[\nabla_{\mathbb{H}^n}u=\frac{\nabla
u}{\lambda},\]

\[\Delta_{\mathbb{H}^n}u=\lambda^{-n}\text{div}(\lambda^{n-2}\nabla u)
;\] where $\nabla$ and $\text{div}$ denote the Euclidean  gradient
and divergence in $\mathbb{R}^n$, respectively.

 The
hyperbolic distance $d_{\mathbb{H}^n}(x,y)$ between
$x,y\in\mathbb{B}^n$ in the Poincare ball model is given by the
formula:
\[d_{\mathbb{H}^n}(x,y)=\text{Arccosh}\Big(1+\frac{2|x-y|^2}{(1-|x|^2)(1-|y|^2)}\Big).\]
From this we immediately obtain for $x\in \mathbb{B}^n$,
\[ \begin{aligned} d:=d_{\mathbb{H}^n}(0,x)&=2\text{Arctanh} |x|\\
&=\log (\frac{1+|x|}{1-|x|})\end{aligned}\] which is the distance
from $x\in \mathbb{B}^n$ to the origin. Moreover, the geodesic
lines passing through the origin are the diameters of $
\mathbb{B}^n$ along with open arcs of circles in $\mathbb{B}^n$
perpendicular to the boundary at $\infty$, $\partial
\mathbb{B}^n=\mathbb{S}^{n-1}= \{x\in \mathbb{R}^n: |x|=1\}$.

The hyperbolic volume element is given by :
\[dV=\lambda^n(x) dx=\Big(\frac{2}{1-r^2}\Big)^nr^{n-1}drd\sigma\]  where $dx$ denotes the  Lebesgue measure in $\mathbb{B}^n$ and $d\sigma$ is the normalized surface measure on $\mathbb{S}^{n-1}$.

A hyperbolic ball in $\mathbb{B}^n$ with center $0$ and hyperbolic
radius $R\in (0, \infty )$ is defined by
\[B_R(0)=\{x\in\mathbb{B}^n \mid d_{\mathbb{H}^n}(0,x)<R\};\]
and note that  $B_R(0)$ is also Euclidean ball with center $0$ and
radius $ S=\text{tanh}\frac{R}{2}\in (0,1)$.

Note that we have the following two relations for the distance
function $d= \log (\frac{1+|x|}{1-|x|})$

\[\begin{aligned}
 |\nabla_{\mathbb{H}^n}d|\, &=1,\\
 \Delta_{\mathbb{H}^n} d\, &\ge
\frac{n-1}{d}, \quad x \neq 0.\end{aligned}\]
\medskip

We are now ready to give a sharp version of Theorem 2.1 above in the
hyperbolic space $\mathbb{H}^n$. Here  $\rho$ is chosen to be the
distance function from the origin in the Poincar\'e ball model for
the hyperbolic space $\mathbb{H}^n$.
\begin{theorem} Let $\phi\in C_0^{\infty}(\mathbb{H}^n)$, $d=\log(\frac{1+|x|}{1-|x|})$, $n\ge 2$, $1<p<\infty$ and
$\alpha>-n$. Then we have:

\begin{equation}\int_{\mathbb{H}^n} d^{\alpha+p}|\nabla_{\mathbb{H}^n} d\cdot \nabla_{\mathbb{H}^n}\phi|^p dV \ge \Big(\frac{n+\alpha}{p}\Big)^p \int_{\mathbb{H}^n}
d^{\alpha}|\phi|^pdV\end{equation} where the constant  $
\big(\frac{n+\alpha}{p}\big)^p$ is sharp.
\end{theorem}
\proof
 The inequality follows from Theorem 2.1. We
 show that
$\big(\frac{n+\alpha}{p}\big)^p$ is the best constant in (2.4):
\[\begin{aligned}C_H:&=\inf_{0\neq\phi\in C_0^{\infty}(\mathbb{H}^n)}\frac{ \int_{\mathbb{H}^n}
d^{\alpha+p}|\nabla_{\mathbb{H}^n}d\cdot
\nabla_{\mathbb{H}^n}\phi|^p dV} {
\int_{\mathbb{H}^n} d^{\alpha}|\phi|^pdV}\\
&=\Big(\frac{n+\alpha}{p}\Big)^p\end{aligned}\]
 It is
clear that

\begin{equation}
\big(\frac{n+\alpha}{p}\big)^p\le\frac{ \int_{\mathbb{H}^n}
d^{\alpha+p}|\nabla_{\mathbb{H}^n}d\cdot
\nabla_{\mathbb{H}^n}\phi|^pdV} { \int_{\mathbb{H}^n}
d^{\alpha}|\phi|^pdV}
\end{equation}
holds for all $\phi\in C_0^{\infty}(\mathbb{H}^n)$. If we pass to
the inf in (2.5) we get that $\big(\frac{n+\alpha}{p}\big)^p\le
C_H$. We only need to show that
$C_H\le\big(\frac{n+\alpha}{p}\big)^p$ and for this we use the
following family of radial functions

\begin{equation} \phi_{\epsilon}(d)=
\begin{cases}
 d^{\frac{n+\alpha}{p}+\epsilon} &\quad\text{if}  \quad d\in [0, 1], \\
d^{-(\frac{n+\alpha}{p}+\epsilon)} &\quad \text{if} \quad d>1,
\end{cases}
\end{equation} where $\epsilon>0$.  Notice that $\phi_{\epsilon}(d)$
can be approximated by smooth functions with compact support in
$\mathbb{H}^n$.

 A direct computation shows that
\[d^{\alpha+p}|\nabla_{\mathbb{H}^n}d\cdot\nabla_{\mathbb{H}^n}\phi_{\epsilon}|^p =
\begin{cases}\big(\frac{n+\alpha}{p}+\epsilon\big)^p d^{n+2\alpha+p\epsilon} &\quad\text{if}  \quad d\in [0, 1], \\
\big(\frac{n+\alpha}{p}+\epsilon\big)^pd^{-n-p\epsilon}
&\quad\text{if} \quad d>1. \end{cases} \]  Let us denote by
$\mathbb{B}_1=\{x\in \mathbb{H}^n: d\le 1\}$ the unit ball with
respect to the distance $d$.  Hence
\[\int_{\mathbb{H}^n} d^{\alpha}
|\phi_{\epsilon}|^pdV=\int_{\mathbb{B}_1}
d^{n+2\alpha+p\epsilon}dV+\int_{\mathbb{H}^n\setminus \mathbb{B}_1}
d^{-n-p\epsilon}dV\]  and then we have
\[\begin{aligned}\big(\frac{n+\alpha}{p}+\epsilon\big)^p\int_{\mathbb{H}^n}d^{\alpha}
|\phi_{\epsilon}|^pdV
&=\big(\frac{n+\alpha}{p}+\epsilon\big)^p\Big[\int_{\mathbb{B}_1}d^{n+2\alpha+p\epsilon}dV+
\int_{\mathbb{H}^n\setminus
\mathbb{B}_1}d^{-n-p\epsilon}dV\Big]\\
&=\int_{\mathbb{H}^n}d^{\alpha+p}|\nabla_{\mathbb{H}^n}d\cdot\nabla_{\mathbb{H}^n}\phi_{\epsilon}|^pdV.\end{aligned}\]
On the other hand

\[\begin{aligned}
\frac{\big(\frac{n+\alpha}{p}+\epsilon\big)^p}{C_H}
\int_{\mathbb{H}^n}d^{\alpha+p}|\nabla_{\mathbb{H}^n}d\cdot\nabla_{\mathbb{H}^n}\phi_{\epsilon}|^pdV
&\ge
\big(\frac{n+\alpha}{p}+\epsilon\big)^p\int_{\mathbb{H}^n}d^{\alpha}
|\phi_{\epsilon}|^pdV
\\ &=\int_{\mathbb{H}^n}d^{\alpha+p}|\nabla_{\mathbb{H}^n}d\cdot\nabla_{\mathbb{H}^n}\phi_{\epsilon}|^pdV.\end{aligned}\]
It is clear that  $\big(\frac{n+\alpha}{p}+\epsilon \big)^p\ge C_H$
and letting $\epsilon \longrightarrow 0$ we obtain
$\big(\frac{n+\alpha}{p} \big)^p\ge C_H$. Therefore
$C_H=\big(\frac{n+\alpha}{p} \big)^p$. \qed
\medskip

We now prove an improved $L^2$ weighted Hardy inequality involving
two weight functions $\rho$ and $\delta$ modeled on distance
functions from a point and distance to the boundary  of a domain
$\Omega$ with smooth boundary.

\begin{theorem} Let $M$ be a complete noncompact Riemannian manifold of dimension $n>1$. Let $\rho$  and $\delta$ be  nonnegative functions on $M$ such
that $|\nabla
 \rho|=1$, $\Delta\rho\ge \frac{C}{\rho}$ and $-\text{div}(\rho^{1-C}\nabla \delta)\ge 0$ in the sense of distribution, where
$C>1$. Then we have:
\begin{equation}\int_{\Omega}\rho^{\alpha}|\nabla\phi|^2dx\ge
\Big(\frac{C+\alpha-1}{2}\Big)^2\int_{\Omega}\rho^{\alpha}\frac{\phi^2}{\rho^2}dx+\frac{1}{4}\int_{\Omega}\rho^{\alpha}\frac{|\nabla
\delta|^2}{\delta^2}\phi^2 dx\end{equation}
 for all $\phi\in C_0^{\infty}(\Omega\setminus
\rho^{-1}\{0\})$, $ \alpha\in \mathbb{R}$ and $C+\alpha-1>0$.

\proof
 Let $\phi\in C_0^{\infty}$ and define $\psi=\rho^{\beta}\phi$  where $\beta<0$. A direct calculation
 shows that
 \begin{equation}
|\nabla\phi|^2=\beta^2\rho^{2\beta-2}|\nabla \rho|^2\psi^2+2\beta
\rho^{2\beta-1}\psi\nabla \rho\cdot
\nabla\psi+\rho^{2\beta}|\nabla\psi|^2.
\end{equation}
Multiplying both sides of (2.8) by the $\rho^{\alpha}$ and applying
integration by parts over $M$ gives

\begin{equation}\begin{aligned}
\int_M \rho^{\alpha}|\nabla\phi|^2dx &=\beta^2\int_M
\rho^{\alpha+2\beta-2}\psi
^2dx-\frac{\beta}{\alpha+2\beta}\int_M
\Delta(\rho^{\alpha+2\beta})\psi^2dx\\&+\int_M \rho^{\alpha+2\beta}|\nabla\psi|^2dx\\
& \ge -\beta^2-\beta(\alpha+C-1)\int_M \rho^{\alpha-2} \phi
^2dx+\int_M \rho^{\alpha+2\beta}|\nabla\psi|^2dx.
\end{aligned}
\end{equation} Choosing
\[\beta=\frac{1-\alpha-C}{2}\]  gives the following
\begin{equation}\int_M \rho^{\alpha}|\nabla\phi|^2dx \ge
\int_M\rho^{\alpha}\frac{\phi^2}{\rho^2}dx+\int_M
\rho^{1-C}|\nabla\psi|^2dx.\end{equation} We now focus on the
second term on the right-hand side of this inequality.  Let us
define a new variable $\varphi (x):=\delta(x)^{-1/2}\psi(x)$ where
$ \delta(x)$ is a nonnegative function and $\delta(x)\in
C_0^2(M)$. It is clear that

\[|\nabla\psi|^2=\frac{1}{4}\frac{\varphi^2}{\delta}|\nabla
\delta|^2+\varphi \nabla\delta\cdot \nabla\varphi+\delta|\nabla
\varphi|^2.\] Therefore

\[
\begin{aligned}\int_M \rho^{1-C}|\nabla\psi|^2dx &\ge
\frac{1}{4}\int_M \rho^{1-C}\frac{\varphi^2}{\delta}|\nabla
\delta|^2+\int_M
\rho^{1-C}\varphi \nabla\delta\cdot \nabla\varphi\\
&= \frac{1}{4}\int_M \rho^{1-C}\frac{|\nabla \delta|^2}{\delta^2}
\psi^2dx-\frac{1}{2}\int_M div (\rho^{1-C}\nabla
\delta)\varphi^2dx.\end{aligned}\] Since $-div(\rho^{1-C}\nabla
\delta)\ge 0$  and $\psi=\rho^{\frac{C+\alpha-1}{2}}\phi$  then we
get
\begin{equation} \int_M \rho^{1-C}|\nabla\psi|^2dx \ge
\frac{1}{4}\int_M \rho^{\alpha}\frac{|\nabla \delta|^2}{\delta^2}
\phi^2dx.\end{equation} Substituting (2.11) into (2.10) gives the
desired inequality: \[\int_M\rho^{\alpha}|\nabla\phi|^2dx\ge
\Big(\frac{C+\alpha-1}{2}\Big)^2\int_M\rho^{\alpha}\frac{\phi^2}{\rho^2}dx+\frac{1}{4}\int_M\rho^{\alpha}\frac{|\nabla
\delta|^2}{\delta^2}\phi^2 dx.\] \qed
\end{theorem}
\medskip

Our next goal is to find model functions which satisfies the
assumption of the above theorem. A straightforward computation
shows that $\delta=\log(\frac{R}{\rho})$ satisfies the
differential inequality $-\text{div}(\rho^{1-C}\nabla \delta)\ge
0$. As a consequence of Theorem 2.3 we have the following
weighted $L^2$-Hardy-type inequality on the hyperbolic space
$\mathbb{H}^n$ which has a logarithmic remainder term. The sharpness of the constant $(\frac{n+\alpha-2}{2})^2$ follows as in \cite{Kombe-Ozaydin} Theorem 3.1.
\medskip

\begin{corollary} Let $\Omega$ be a bounded domain with smooth boundary $\partial\Omega$ in $\mathbb{H}^n$.
 Let $\rho=d=\log(\frac{1+|x|}{1-|x|})$ and $\delta:=\log(\frac{R}{d})$, $R> \sup_{\Omega}\big(d\big)$, $\alpha\in \mathbb{R}$,
$n+\alpha-2>0$. Then we have:
\begin{equation}\int_{\Omega} d^{\alpha}|\nabla_{\mathbb{H}^n} \phi|^2dV
\ge \Big(\frac{n+\alpha-2}{2}\Big)^2 \int_{\Omega}
d^{\alpha}\frac{\phi^2}{d^2}dV+\frac{1}{4}\int_{\Omega}d^{\alpha}\frac{\phi^2}{d^2(\log\frac{R}{d})^2}dV\end{equation}
for all  $\phi\in C_0^{\infty}(\Omega)$ and the constant
$\big(\frac{n+\alpha-2}{2}\big)^2$ is sharp.
\end{corollary}
\medskip

Let $B_R=\{x\in \mathbb{B}^n\mid d<R\}$ be a hyperbolic ball with
center $0$ and hyperbolic radius $R$. It is clear that $\delta:=R-d$
is the distance function of the point $x\in B_R$ to the boundary of
$B_R$ and satisfies the differential inequality in Theorem 2.3.
Therefore we have:
\begin{corollary} Let $B_R$ be a hyperbolic ball with center $0$ and hyperbolic radius $R$.
 Let $d=\log(\frac{1+|x|}{1-|x|})$ and $\delta:=R-d$, $\alpha\in \mathbb{R}$,
$n+\alpha-2>0$. Then we have:
\begin{equation}\int_{B_R} d^{\alpha}|\nabla_{\mathbb{H}^n} \phi|^2dV
\ge \Big(\frac{n+\alpha-2}{2}\Big)^2 \int_{B_R}
d^{\alpha}\frac{\phi^2}{d^2}dV+\frac{1}{4}\int_{B_R}d^{\alpha}\frac{\phi^2}{(R-d)^2}dV\end{equation}
for all  $\phi\in C_0^{\infty}(B_R)$ and the constant
$\big(\frac{n+\alpha-2}{2}\big)^2$ is sharp.
\end{corollary}
\medskip

\noindent{\textbf{Hardy-Sobolev-Poincar\'e inequalities}}. The
following sharp form of the Sobolev inequality on the hyperbolic
space $\mathbb{H}^n$ is due to \cite{Hebey}. It states that for all
$\phi\in C_0^{\infty}(\mathbb{H}^n)$:

\begin{equation}\int_{
\mathbb{H}^n }|\nabla_{\mathbb{H}^n}\phi|^2dV\ge
\frac{n(n-2)}{4}|\mathbb{S}^n|^{\frac{2}{n}}\Big(\int_ {\mathbb
{H}^n}|\phi|^{\frac{2n}{n-2}}dV\Big)^{\frac{n-2}{n}}+\frac{n(n-2)}{4}\int_{
\mathbb{H}^n} \phi^2dV
 \end{equation}where $
 \phi\in C_0^{\infty}(\mathbb{H}^n)$. Here $A_n=\frac{n(n-2)}{4}|\mathbb{S}|^{\frac{2}{n}}$ is the sharp
constant for the Sobolev inequality on $\mathbb{R}^n$,
$|\mathbb{S}^n|$ is the volume of the $n$-dimensional unit sphere in
$\mathbb{R}^{n+1}$ and the constant $ B_n=\frac{n(n-2)}{4}$ is sharp
for $n\ge 4$. Recently, sharp form of the inequality (2.13) in
three dimensional hyperbolic space $\mathbb{H}^n$ has been proved by
Benguria, Frank and Loss \cite{Benguria}.
\medskip

The Sobolev inequality (2.14) and  Hardy inequality
\cite{Kombe-Ozaydin} yield the following Hardy-Sobolev inequality in
$\mathbb{H}^n$.

\begin{corollary} Let $\phi\in C_0^{\infty}(\mathbb{H}^n)$, $d=\log(\frac{1+|x|}{1-|x|})$ and $n\ge 3$. Then we have:
\[\int_{ \mathbb{H}^n }|\nabla_{\mathbb{H}^n}\phi|^2dV\ge
\Big(\frac{n-2}{2}\Big)^{\frac{2s}{p^{*}(s)}}
\Big(\frac{n(n-2)}{4}|\mathbb{S}^n|^{\frac{2}{n}}\Big)
  ^{\frac{n(2-s)}{2(n-s)}}\Big(\int_ {\mathbb
{H}^n}\frac{|\phi|^{p^{*}(s)}}{d^s}dV\Big)^{\frac{2}{p^{*}(s)}}
\] where $0\le s\le 2$ and  $p^{*}(s)=2 (\frac{n-s}{n-2})$.
\end{corollary}
\medskip

Before we state and prove our next theorem, we first recall the
(Euclidean) weighted Sobolev inequality of Fabes-Kenig-Serapino
\cite{Kenig} which plays an important role in our proof. They proved
the following inequality :
\begin{equation}\Big(\frac{1}{w(B_r)}\int_{B_r}|\nabla\phi|^p
w(x)dx\Big)^{1/p}\ge \frac{1}{c(\text{diam}
B_r)}\Big(\frac{1}{w(B_r)}\int_{B_r}|\phi|^{kp}w(x)dx\Big)^{1/kp}
\end{equation} where $B_r$ is a ball in
$\mathbb{R}^n$,  $\phi\in C_0^{\infty}(B_r)$,
$w(B_r)=\int_{B_r}w(x)dx$, $1<p<\infty$, $1\le k\le
\frac{n}{n-1}+\epsilon$, $\epsilon>0$ and the weight function $w$
belongs to Muckenhoupt's class $A_p$. In particular, if the weight
function $w$ belongs to Muckenhoupt's class $A_2$ then $k$ can be
taken equal to $\frac{n}{n-1}+\epsilon$ and this is sharp. Recall
that a weight function $w$ belongs to Muckenhoupt's class $A_p$
($1<p<\infty$) if
\[\text{sup}\Big(\frac{1}{|B|}\int_{B}w(x)dx\Big)\Big(\frac{1}{|B|}\int_{B}w(x)^{\frac{1}{1-p}}dx\Big)^{p-1}=
C_{p,w}<\infty, \]  where the supremum is taken over all balls $B$
in $\mathbb{R}^n$ (see \cite{Stein}).
\medskip

Motivated by the classical work of Brezis and V\'azquez \cite{Brezis},
our next theorem shows that sharp weighted Hardy inequality on the
hyperbolic space $\mathbb{H}^n$ can be improved by a weighted
Sobolev term.

\begin{theorem}
Let  $\phi\in C_0^{\infty}(\mathbb{H}^n)$, $d=\log
(\frac{1+|x|}{1-|x|})$, $\alpha\in\mathbb{R}$, $n>2$ and
$n+\alpha-2>0$. Then we have : \begin{equation}\int_{ \mathbb{H}^n
}d^{\alpha}|\nabla_{\mathbb{H}^n}\phi|^2dV\ge
\Big(\frac{n+\alpha-2}{2}\Big)^2\int_ {\mathbb
{H}^n}d^{\alpha}\frac{\phi^2}{d^2}dV+ \tilde{c}\Big(\int_{
\mathbb{H}^n} d^{\frac{(2-n)(2-q)+\alpha q}{2}}\phi^qdx\Big)^{2/q}
 \end{equation} where $2\le q\le \frac{2n}{n-1}+2\epsilon$, $\epsilon>0$, $\tilde{c}=
\frac{2^{n-2}}{c^2}\big(\frac{|\mathbb{S}^n|}{2}\big)^{\frac{q-2}{q}}
$, $c>0$ and  the constant $(\frac{n+\alpha-2}{2})^2$ is sharp.
\end{theorem}

\proof  Let $\phi\in C_0^{\infty}$ and define $\psi=d^{-\beta}\phi$
where $\beta<0$. A direct calculation shows that
\begin{equation}
d^{\alpha}|\nabla\phi|^2\lambda^{n-2}
=\beta^2d^{\alpha+2\beta-2}|\nabla d|^2\psi^2\lambda^{n-2}+ 2\beta
d^{\alpha+2\beta-1}\psi \lambda^{n-2}\nabla
d\cdot\nabla\psi+d^{\alpha+2\beta}|\nabla\psi|^2 \lambda^{n-2}.
\end{equation}
It is easy to see that
\[|\nabla d|^2=\lambda^2\]
and integrating (2.17) over $ \mathbb{B}^n$, we get
\begin{equation}
\begin{aligned}\int_{\mathbb{B}^n}d^{\alpha}|\nabla\phi|^2 \lambda^{n-2}dx&=\int_{
\mathbb{B}^n}\beta^2 d^{\alpha+2\beta-2}\psi^2
\lambda^ndx+\int_{\mathbb{B}^n} 2\beta d^{\alpha+2\beta-1}\psi
\lambda^{n-2}\nabla d\cdot\nabla\psi dx\\&+\int_{\mathbb{B}^n}
d^{\alpha+2\beta}|\nabla\psi|^2 \lambda^{n-2}dx.\\
\end{aligned}
\end{equation}
Applying integration by parts to the middle integral on the
right-hand side of (2.18), we obtain

\begin{equation}
\begin{aligned}
\int_{\mathbb{B}^n}d^{\alpha}|\nabla\phi|^2
\lambda^{n-2}dx&=\int_{ \mathbb{B}^n}\beta^2
d^{\alpha+2\beta-2}\psi^2 \lambda^ndx
 -\frac{\beta}{\alpha+2\beta} \int_{\mathbb{B}^n}\text{div}\big(\lambda^{n-2}\nabla (d^{2\beta+\alpha})\big)dx \\
 &+\int_{\mathbb{B}^n}
d^{\alpha+2\beta}|\nabla\psi|^2 \lambda^{n-2}dx.
\end{aligned}
\end{equation}
One can show that
\begin{equation}
\begin{aligned}
&-\frac{\beta}{\alpha+2\beta}
\int_{\mathbb{B}^n}\text{div}\big(\lambda^{n-2}\nabla
(d^{2\beta+\alpha})\big)dx \\ =& -\beta (2\beta
+\alpha-1)\int_{\mathbb{B}^n}d^{2\beta+\alpha-2}\lambda^n\psi^2dx-\beta
\int_{\mathbb{B}^n}d^{2\beta+\alpha-1}\lambda^{n-2}\psi^2 (\Delta d)dx\\
&-\beta
(n-2)\int_{\mathbb{B}^n}d^{2\beta+\alpha-1}\lambda^{n-3}(\nabla
d\cdot \nabla \lambda) dx.
\end{aligned}
\end{equation}
A direct computation shows that
\[\Delta d=\lambda^2r+\frac{n-1}{r}\lambda\] and

\[\nabla d\cdot \nabla \lambda=\lambda^3r.\] Substituting these above
\begin{equation}
\begin{aligned}
&-\frac{\beta}{\alpha+2\beta}
\int_{\mathbb{B}^n}\text{div}\big(\lambda^{n-2}\nabla
(d^{2\beta+\alpha})\big)dx \\ =& -\beta (2\beta
+\alpha-1)\int_{\mathbb{B}^n}d^{2\beta+\alpha-2}\lambda^n\psi^2dx-(2\beta
+\alpha)\int_{\mathbb{B}^n}
d^{2\beta+\alpha-1}\lambda^n\big(\frac{(n-1)(\lambda
r^2+1)}{\lambda r}\big)\psi^2 dx.
\end{aligned}
\end{equation}
We can easily show that \[\frac{\lambda r^2+1}{\lambda r}\ge
\frac{1}{d}.\] If $2\beta+\alpha<0$ then we have

\begin{equation}
-\frac{\beta}{\alpha+2\beta}
\int_{\mathbb{B}^n}\text{div}\big(\lambda^{n-2}\nabla
(d^{2\beta+\alpha})\big)dx \ge
-\beta(2\beta+\alpha+n-2)\int_{\mathbb{B}^n}
d^{2\beta+\alpha-2}\lambda^n\psi^2dx.
\end{equation}
Now we substitute (2.22) into (2.19) and we get

\[\int_{
\mathbb{B}^n}d^{\alpha}|\nabla \phi|^2\lambda^{n-2}\ge
(-\beta^2-\beta(\alpha+n-2)) \int_{
\mathbb{B}^n}d^{2\beta+\alpha-2}\psi^2 \lambda^n
dx+\int_{\mathbb{B}^n} d^{\alpha+2\beta}|\nabla\psi|^2
\lambda^{n-2}dx.\]

Note that the function $\beta\longrightarrow
-\beta^2-\beta(\alpha+n-2)$ attains the maximum for
$\beta=\frac{2-\alpha-n}{2}$, and this maximum is equal to
$(\frac{n+\alpha-2}{2})^2$. Therefore we have the following
inequality

\[\int_{
\mathbb{B}^n}d^{\alpha}|\nabla \phi|^2\lambda^{n-2}dx\ge
\Big(\frac{n+\alpha-2}{2}\Big)^2\int_
{\mathbb{B}^n}d^{\alpha}\frac{\phi^2}{d^2} \lambda^n
dx+\int_{\mathbb{B}^n} d^{2-n}|\nabla\psi|^2 \lambda^{n-2}dx.\] Using the fact $d\le \lambda r$ we get\\
\begin{equation}\int_{ \mathbb{B}^n}d^{\alpha}|\nabla
\phi|^2\lambda^{n-2}dx\ge\Big(\frac{n+\alpha-2}{2}\Big)^2\int_
{\mathbb{B}^n}d^{\alpha}\frac{\phi^2}{d^2} \lambda^n
dx+\int_{\mathbb{B}^n} r^{2-n}|\nabla\psi|^2 dx.\end{equation}
 Notice that the weight
function $r^{2-n}$ is in the Muckenhoupt $A_2$ class. We now apply
weighted Sobolev inequality (2.15) to the second integral term on
the right hand side of (2.23) and obtain

\[\begin{aligned}
\int_{\mathbb{B}^n}d^{\alpha}|\nabla\phi|^2\lambda^{n-2}dx
&\ge\Big(\frac{n+\alpha-2}{2}\Big)^2\int_{\mathbb{B}^n}
d^{\alpha}\frac{\phi^2}{d^2}\lambda^ndx+ c_1\Big(\int_{
\mathbb{B}^n} r^{2-n}\psi^qdx\Big)^{2/q}\\
&\ge \Big(\frac{n+\alpha-2}{2}\Big)^2\int_{\mathbb{B}^n}
d^{\alpha}\frac{\phi^2}{d^2}\lambda^ndx+c_1\Big(\int_{ \mathbb{B}^n}
r^{2-n}d^{\frac{(n+\alpha-2)q}{2}}\phi^qdx\Big)^{2/q}\end{aligned}\]
where $q>2$ and
$c_1=\frac{1}{c^2}\big(\frac{|\mathbb{S}^n|}{2}\big)^{1-\frac{1}{k}}$.
Furthermore,  using the inequality $2r\le d\le \lambda r$, we get
\[\int_{\mathbb{B}^n}d^{\alpha}|\nabla\phi|^2\lambda^{n-2}dx\ge\Big(\frac{n+\alpha-2}{2}\Big)^2\int_{\mathbb{B}^n}
d^{\alpha}\frac{\phi^2}{d^2}\lambda^ndx+ \tilde{c}\Big(\int_{
\mathbb{B}^n} d^{\frac{(n-2)(q-2)+\alpha q}{2}}\phi^qdx\Big)^{2/q}\]
where $q>2$ and $\tilde{c}=
\frac{2^{n-2}}{c^2}\big(\frac{|\mathbb{S}^n|}{2}\big)^{\frac{q-2}{q}}
$.  This completes the proof.\qed

\section{Rellich-type inequalities}
In this section we prove weighted Rellich-type inequality and its
improved versions which connects first to the second order
derivatives. The following is the weighted analogue of (1.4) in the
setting of Riemannian manifold $M$.
\begin{theorem} Let $M$ be a complete  Riemannian manifold of dimension
$n>1$. Let $\rho$ be a nonnegative function on $M$ such that
$|\nabla \rho|=1$ and  $\Delta\rho\ge \frac{C}{\rho}$ in the sense
of distribution where $ C>1$. Then the following inequality is
valid:
\begin{equation}
\int_{M}\rho^{\alpha}|\Delta\phi|^2dV \ge \frac{(C+1-\alpha)^2}{4}
\int_{M}\rho^{\alpha}\frac{|\nabla\phi|^2}{\rho^{2}}dV.
\end{equation}
for all compactly supported smooth function $\phi\in
C_0^{\infty}(M\setminus\rho^{-1}\{0\})$, $\frac{7-C}{3}<\alpha<2$
\end{theorem}
\proof A straightforward computation shows that
\begin{equation}\Delta \rho^{\alpha-2}\le(\alpha-2)(C+\alpha-3)\rho^{\alpha-4}.\end{equation}
Multiplying both sides of (3.2) by $\phi^2$ and integrating over
$M$,  we obtain
\begin{equation}\begin{aligned}(C+\alpha-3)(\alpha-2)\int_M \rho^{\alpha-4}\phi^2dV &\ge \int_M \rho^{\alpha-2}
\Delta(\phi^2)dV\\& = \int_M
\rho^{\alpha-2}(2|\nabla\phi|^2+2\phi\Delta\phi)dV.
\end{aligned}
\end{equation}
Therefore
\begin{equation}-\int_M(\phi\Delta\phi)\rho^{\alpha-2}\ge \int_M
\rho^{\alpha-2}|\nabla\phi|^2dV-\frac{(C+\alpha-3)(\alpha-2)}{2}\int_M
\rho^{\alpha-4}\phi^2dV.\end{equation}

Let us apply Young's inequality to expression $-\int_M
\rho^{\alpha-2}\phi\Delta\phi\,dx$
\begin{equation}-\int_M \rho^{\alpha-2}\phi\Delta\phi dV
\le \epsilon \int_M
\rho^{\alpha-4}\phi^2dV+\frac{1}{4\epsilon}\int_M
\rho^{\alpha}|\Delta\phi|^2dV\end{equation}
 where $\epsilon>0$ and will be chosen later. Combining (3.5) and
 (3.4) we get

 \begin{equation}
\int_M \rho^{\alpha-2}|\nabla \phi|^2dV\le
\Big(\epsilon+\frac{(C+\alpha-3)(\alpha-2)}{2}\Big)\int_M
\rho^{\alpha-4}\phi^2dV+\frac{1}{4\epsilon} \int_M
\rho^{\alpha}|\Delta\phi|^2dV.
\end{equation}
Notice that the case of
$\epsilon+\frac{(C+\alpha-3)(\alpha-2)}{2}<0$ gives the Rellich
inequality (1.6). Therefore we only need to consider the cases:
$\epsilon+\frac{(C+\alpha-3)(\alpha-2)}{2}=0$  and
$\epsilon+\frac{(C+\alpha-3)(\alpha-2)}{2}>0$, respectively. The
first case gives the following inequality:

\begin{equation}
\int_M \rho^{\alpha}|\Delta\phi|^2dV\ge
2(C+\alpha-3)(2-\alpha)\int_M \rho^{\alpha-2}|\nabla \phi|^2dV.
\end{equation} If $\epsilon+\frac{(C+\alpha-3)(\alpha-2)}{2}>0$ then
we apply the Rellich inequality (1.6) to the first term on the
right hand side of (3.6) and get

\begin{equation}
\int_M \rho^{\alpha-2}|\nabla \phi|^2dV\le P_{C,\alpha}(\epsilon)
\int_M \rho^{\alpha}|\Delta\phi|^2dV
\end{equation}
where
\[P_{C,\alpha}(\epsilon)=\frac{16\epsilon}{(C+\alpha-3)^2(C-\alpha+1)^2}+\frac{8(\alpha-2)}{(C+\alpha-3)(C-\alpha+1)^2}+\frac{1}{4\epsilon}.\]  Note that the function $P_{C,\alpha}(\epsilon)$  attains the minumum for $\epsilon=\frac{(C+\alpha-3)(C-\alpha+1)}{8}$, and this minimum is equal to $\frac{4}{(C-\alpha+1)^2}$. Therefore we have the following inequality:

\[\int_M
\rho^{\alpha}|\Delta\phi|^2dV \ge \frac{(C-\alpha+1)^2}{4}\int_M
\rho^{\alpha-2}|\nabla \phi|^2dV.\] \qed
\medskip

We are now ready to give a sharp version of Theorem 3.1 above in the
hyperbolic space $\mathbb{H}^n$. Here  $\rho$ is chosen to be the
distance function from the origin in the Poincar\'e ball model for
the hyperbolic space.
\begin{theorem} Let $\phi\in C_0^{\infty}(\mathbb{H}^n)$, $d=\log(\frac{1+|x|}{1-|x|})$, $n> 2$,
$\frac{8-n}{3}<\alpha<2$. Then we
have:\begin{equation}\int_{\mathbb{H}^n}
d^{\alpha}|\Delta_{\mathbb{H}^n}\phi|^2dV \ge
\frac{(n-\alpha)^2}{4}\int_{\mathbb{H}^n}
d^{\alpha}\frac{|\nabla_{\mathbb{H}^n}
\phi|^2}{d^2}dV.\end{equation} \proof The inequality follows from
Theorem 3.1. To show that the constant $(\frac{n-\alpha}{2})^2$ is
sharp, we use the following family of functions as in
\cite{Kombe-Ozaydin}:

\[\phi_{\epsilon}(d)=\begin{cases}
 -(\frac{n+\alpha-4}{2}+\epsilon)\big(d-1\big)+1 &\quad\text{if}  \quad d\in [0,1],\\
d^{-(\frac{n+\alpha-4}{2}+\epsilon)} &\quad \text{if} \quad d>1.
\end{cases}
\]
Notice that $\phi_{\epsilon}(d)$ can be well approximated by smooth
functions with compact support in $\mathbb{H}^n$ and direct
computation shows that $\frac{(n-\alpha)^2}{4}$ is the best constant
in (3.9):
\[\begin{aligned}\frac{\big(n-\alpha\big)^2}{4}=\lim_{\epsilon\longrightarrow 0}\frac{ \int_{\mathbb{H}^n}
d^{\alpha}|\Delta_{\mathbb{H}^n}\phi_{\epsilon}|^2 dV} {
\int_{\mathbb{H}^n}
d^{\alpha}\frac{|\nabla_{\mathbb{H}^n}\phi_{\epsilon}|^2}{d^2}dV}.
\end{aligned}\]
\end{theorem}
\medskip

 The following inequality is an improved
version of the Rellich-type inequality (3.1) for bounded domains.
\begin{theorem}Let $\Omega$ be a bounded domain with smooth boundary $\partial\Omega$ in a  complete  Riemannian manifold of dimension $n>1$.
  Let $\rho$ be a nonnegative function on $M$ such that
  $|\nabla\rho|=1$,
 $\Delta\rho\ge \frac{C}{\rho}$ and $-div(\rho^{1-C}\nabla\delta)\ge 0$ in the sense of distribution, where $ C>1$.  Then the following inequality is
valid:
\begin{equation}
\int_{\Omega}\rho^{\alpha}|\Delta\phi|^2dV \ge
\frac{(C+1-\alpha)^2}{4}
\int_{\Omega}\rho^{\alpha}\frac{|\nabla\phi|^2}{\rho^{2}}dV+K(C,
\alpha)\int_{\Omega}\rho^{\alpha-2}\frac{|\nabla
\delta|^2}{\delta^2}\phi^2dV
\end{equation}
for all compactly supported smooth function $\phi\in
C_0^{\infty}(M\setminus\rho^{-1}\{0\})$, $\frac{7-C}{3}<\alpha<2$
and $K(C, \alpha)=\frac{(C+1-\alpha)(C+3\alpha-7)}{16}$. \proof The
proof is similar to the proof of Theorem 3.1. The only difference is
that we apply improved Hardy-type inequality (2.7) to the first term
on the right hand side of (3.6).\qed
\end{theorem}
\medskip

The following corollaries are the direct consequences of the Theorem
3.3.

\begin{corollary} Let $\Omega$ be a bounded domain with smooth boundary in $\mathbb{H}^n$.
 Let $d=\log(\frac{1+|x|}{1-|x|})$ and $\delta:=\log(\frac{R}{d})$ and $R> \sup_{\Omega}\big(d\big)$. Then the following inequality is valid:
\begin{equation}
\int_{\Omega}d^{\alpha}|\Delta_{\mathbb{H}^n}\phi|^2dV \ge
\frac{(n-\alpha)^2}{4}
\int_{\Omega}d^{\alpha}\frac{|\nabla_{\mathbb{H}^n}\phi|^2}{d^{2}}dV+K(C,
\alpha)\int_{\Omega}d^{\alpha-4}\frac{\phi^2}{(\ln\frac{R}{d})^2}dV
\end{equation}
for all compactly supported smooth function $\phi\in
C_0^{\infty}(\Omega)$, $\frac{8-n}{3}<\alpha<2$ and
$K=\frac{(n-\alpha)(n+3\alpha-8)}{16}$.
\end{corollary}
\medskip

\begin{corollary} Let $B_R$ be a hyperbolic ball with center $0$ and hyperbolic radius $R$.
 Let $d=\log(\frac{1+|x|}{1-|x|})$ and $\delta:=R-d$. Then the following inequality is valid:
\begin{equation}
\int_{\Omega}d^{\alpha}|\Delta_{\mathbb{H}^n}\phi|^2dV \ge
\frac{(n-\alpha)^2}{4}
\int_{\Omega}d^{\alpha}\frac{|\nabla_{\mathbb{H}^n}\phi|^2}{d^{2}}dV+K(C,
\alpha)\int_{\Omega}d^{\alpha-2}\frac{\phi^2}{(R-d)^2}dV
\end{equation}
for all compactly supported smooth function $\phi\in
C_0^{\infty}(\Omega)$, $\frac{8-n}{3}<\alpha<2$ and
$K=\frac{(n-\alpha)(n+3\alpha-8)}{16}$.
\end{corollary}

Using the same argument as in the proof of Theorem 3.1 and improved
Hardy-Sobolev type inequality (2.16), we obtain the following
improved Rellich-Sobolev type inequality on the hyperbolic space
$\mathbb{H}^n$.
\begin{corollary} Let  $\phi\in C_0^{\infty}(\mathbb{H}^n)$ and $d=\log
(\frac{1+|x|}{1-|x|})$. Then the following inequality is valid:
\begin{equation}
\int_{\mathbb{H}^n}d^{\alpha}|\Delta_{\mathbb{H}^n}\phi|^2dV \ge
\frac{(n-\alpha)^2}{4}
\int_{\mathbb{H}^n}d^{\alpha}\frac{|\nabla_{\mathbb{H}^n}\phi|^2}{d^{2}}dV+K\Big(\int_{
\mathbb{H}^n} d^{\frac{(n-2)(q-2)+(\alpha-2)
q}{2}}\phi^qdx\Big)^{2/q}
\end{equation}
where  $\frac{8-n}{3}<\alpha<2$,
$K=\frac{(n-\alpha)(n+3\alpha-8)2^{n-2}}{4c^2}\big(\frac{|\mathbb{S}^n|}{2}\big)^{\frac{q-2}{q}}
$, $2\le q\le \frac{2n}{n-1}+2\epsilon$, $\epsilon>0$, and $c>0$.
\end{corollary}

\section{Uncertainty Principle Inequality.}
The first and most famous uncertainty principle goes back to
Heisenberg's seminal work  which was developed in the context of
quantum mechanics \cite {Heisenberg}. It says that the position
and momentum of a particle cannot be determined exactly at the
same time but only with an ``uncertainty".  The mathematical
version of this principle (stating that a function and its Fourier
transform can not be well localized simultaneously) was formulated
afterwards by Pauli and Weyl \cite{Weyl} and is sometimes referred
to as the Heisenberg-Pauli-Weyl inequality. Uncertainty principle
type inequalities are  central to harmonic analysis and such
considerations of  the time-frequency domain are crucial in signal
and image processing \cite{Folland-Sitaram}.

In a previous work \cite{Kombe-Ozaydin},  we obtained a
Heisenberg-Pauli-Weyl inequality on a compete non-compact Riemannian
manifold $M$ and found an explicit constant when $M$ is
 the hyperbolic space $\mathbb{H}^n$. In this present paper we
first prove a Heisenberg-Pauli-Weyl inequality for general
Riemannian manifolds which has a better constant than those of
\cite{Kombe-Ozaydin} and then obtain the sharp constant in the
Hyperbolic case. The following is the first result of this section.

\begin{theorem}  Let $M$ be a complete  Riemannian manifold of dimension $n\ge2$. Let $\rho$ be a nonnegative function on $M$ such that
$|\nabla \rho|=1$ and $\Delta\rho\ge \frac{C}{\rho}$ in the sense
of distribution where $C>0$. Then the following inequality holds:
\begin{equation}
\Big(\int_M \rho^2\phi^2dV\Big)\Big(\int_M|\nabla \phi|^2dV\Big)\ge
\frac{(C+1)^2}{4}\Big(\int_M\phi^2 dV\Big)^2
\end{equation}
for all compactly supported smooth function $\phi\in
C_0^{\infty}(M)$. \proof Using the  assumptions $|\nabla \rho|=1$
and $\Delta\rho\ge \frac{C}{\rho}$, we get
\begin{equation}\int_M(\Delta\rho^2)\phi^2dV\ge
(2C+2)\int_M\phi^2dV. \end{equation} By integration by parts and the
Cauchy-Schwarz inequality, we have
\[\Big(\int_M
\rho^2\phi^2dV\Big)\Big(\int_M|\nabla \phi|^2dV\Big)\ge
\frac{(C+1)^2}{4}\Big(\int_M\phi^2 dV\Big)^2.\] This completes the
proof.\qed
\end{theorem}
We now proof a sharp analogue of the Heisenberg-Pauli-Weyl
inequality (1.3) on the hyperbolic space $\mathbb{H}^n$.
\begin{theorem} Let
$\phi\in C_0^{\infty}(\mathbb{H}^n)$, $d=\log(\frac{1+|x|}{1-|x|})$
and $n> 2$. Then
\begin{equation}
\Big(\int_{\mathbb{H}^n} d^2 \phi^2
dV\Big)\Big(\int_{\mathbb{H}^n}|\nabla_{\mathbb{H}^n}
\phi|^2dV\Big)\ge \frac{n^2}{4}\Big(\int_{\mathbb{H}^n}\phi^2
dV\Big)^2.
\end{equation}
Moreover, equality holds in (4.3) if $\phi(x)=Ae^{-\alpha d^2}$
where $A\in \mathbb{R}$, $\alpha=(\frac{n-1}{n-2})\big(n-1+2\pi
\frac{C_{n-2}}{C_n}\big)$ and $C_n=\int_{\mathbb{H}^n}e^{-\alpha
d^2}dV$ . \proof The inequality follows from Theorem 4.1. In order
to achieve equality, inspired by the Euclidean case, we consider
hyperbolic analogues of Gaussians: $\phi(x)=Ae^{-\alpha d^2}$ where
$A\in \mathbb{R}$ and $\alpha>0$. A straightforward but tedious
calculation shows that $\phi(x)=Ae^{-\alpha d^2}$ is the minimizer
where $\alpha=(\frac{n-1}{n-2})\big(n-1+2\pi
\frac{C_{n-2}}{C_n}\big)$ and  $C_n=\int_{\mathbb{H}^n}e^{-\alpha
d^2}dV$.\qed
\end{theorem}
\noindent{\bf Remark}. Note that even though $\phi(x)=Ae^{-\alpha d^2}$ does not have a compact support, it can be approximated by such functions yielding that (4.3) is sharp.
\medskip

There is a natural link between Hardy, Heisenberg-Pauli-Weyl and
Rellich type inequalities. For instance, using the
\textit{Rellich-type inequality II} (3.1) we have the following
second order Heisenberg-Pauli-Weyl inequality.

\begin{theorem} Let $M$ be a complete  Riemannian manifold of dimension $n\ge2$. Let $\rho$ be a nonnegative function on $M$ such that
$|\nabla \rho|=1$ and $\Delta\rho\ge \frac{C}{\rho}$ in the sense of
distribution where $C>1$. Then the following inequality holds:
\begin{equation}
\Big(\int_M \rho^4\phi^2dV\Big)\Big(\int_M|\Delta
\phi|^2dV\Big)\ge\frac{(C+1)^4}{16} \Big(\int_M\phi^2 dV\Big)^2.
\end{equation} for all compactly supported smooth function $\phi\in
C_0^{\infty}(M-\rho^{-1}\{0\})$.
 \proof By the equation (4.2) and Cauchy-Schwarz
inequality, we get
\[\Big(\int_M\rho^4\phi^2dV\Big)^{1/2}\Big(\int_M
\frac{|\nabla\phi|^2}{\rho^2}dV\Big)^{1/2}\ge \frac{C+1}{4}
\int_M\phi^2dV.\] Using the  \textit{Rellich-type inequality II}
(3.1) we obtain the desired inequality:

\[\Big(\int_M\rho^4\phi^2dV\Big)\Big(\int_M|\Delta
\phi|^2dV\Big)\ge\frac{(C+1)^4}{16} \Big(\int_M\phi^2 dV\Big)^2.\]
\end{theorem}
\qed

As an immediate consequence of the Theorem 4.3 we have the following
second order Heisenberg-Pauli-Weyl inequality with an explicit
constant on the hyperbolic space $\mathbb{H}^n$.
\begin{corollary}  Let
$\phi\in C_0^{\infty}(\mathbb{H}^n-\{0\})$,
$d=\log(\frac{1+|x|}{1-|x|})$, $n> 2$ and $\frac{8-n}{3}<\alpha<2$.
Then the following inequality holds:
\begin{equation}
\Big(\int_{\mathbb{H}^n}
d^4\phi^2dV\Big)\Big(\int_{\mathbb{H}^n}|\Delta_{\mathbb{H}^n}
\phi|^2dV\Big)\ge\frac{n^4}{16} \Big(\int_{\mathbb{H}^n}\phi^2
dV\Big)^2.
\end{equation}
\end{corollary}

\bibliographystyle{amsalpha}

\end{document}